\documentclass[12pt,a4paper,reqno]{amsart}
\usepackage{amssymb}
\usepackage{amscd}

\numberwithin{equation}{section}

\theoremstyle{plain}

\newtheorem{theorem}[subsection]{Theorem}
\newtheorem{proposition}[subsection]{Proposition}
\newtheorem{lemma}[subsection]{Lemma}

\newtheorem{conjecture}[subsection]{Conjecture}

\theoremstyle{definition}

\newtheorem{problem}[subsection]{Problem}

\renewcommand{\leq}{\leqslant}
\renewcommand{\geq}{\geqslant}

\newsavebox{\proofbox}
\savebox{\proofbox}{\begin{picture}(7,7)%
  \put(0,0){\framebox(7,7){}}\end{picture}}

\newcommand{\md}[1]{\ensuremath{(\operatorname{mod}\, #1)}}

\newcommand{\mdlem}[1]{\ensuremath{(\mbox{\textup{mod}}\, #1)}}

\newcommand\Z{\mathbb{Z}}
\newcommand\R{\mathbb{R}}

\newcommand\eps{\varepsilon}

\def\proof{\textit{Proof. }}

\def\endproof{\hfill{\usebox{\proofbox}}}

\begin{document}

\title{On a variant of the large sieve}

\author{Ben Green}
\address{Centre for Mathematical Sciences\\
Wilberforce Road\\
     Cambridge CB3 0WA\\
     England
}
\email{b.j.green@dpmms.cam.ac.uk}

\thanks{The author holds a Leverhulme Prize and is grateful to the Leverhulme Trust for their support.}

\begin{abstract}
We introduce a variant of the large sieve and give an example of its use in a sieving problem. Take the interval $[N] = \{1,\dots,N\}$ and, for each odd prime $p \leq \sqrt{N}$, remove or ``sieve out'' by all $n$ whose reduction $n\md{p}$ lies in some interval $I_p \subseteq \Z/p\Z$ of length $(p-1)/2$. Let $A$ be the set that remains: then $|A| \ll_{\eps} N^{1/3 + o(1)}$, a bound which improves slightly on the bound of $|A| \ll N^{1/2}$ which results from applying the large sieve in its usual form. This is a very, very weak result in the direction of a question of Helfgott and Venkatesh, who suggested that nothing like equality can occur in applications of the large sieve unless the unsieved set is essentially the set of values of a polynomial (e.g. $A$ is the set of squares).

Assuming the ``exponent pairs conjecture'' (which is deep, as it implies a host of classical questions including the Lindel\"of hypothesis, Gauss circle problem and Dirichlet divisor problem) we can improve the bound to $|A| \ll N^{o(1)}$. This raises the worry that even reasonably simple sieve problems are connected to issues of which we have little understanding at the present time.
\end{abstract}

\maketitle

\section{Introduction}

The large sieve is, in its purest form, the following analytic inequality.

\begin{theorem}[Large sieve]\label{large-sieve}
Suppose that $N \geq 1$ is an integer and that $\delta \in (0,1)$. Suppose that the points $\theta_1,\dots,\theta_k \in \R/\Z$ are $\delta$-separated, that is to say \begin{equation}\label{separation} |\theta_i - \theta_j | \geq \delta \end{equation} whenever $1 \leq i < j \leq k$. Let $(a_n)_{n \in [N]}$ be any sequence of complex numbers. Then
\begin{equation}\label{l-sieve}\sum_{i=1}^k |\sum_{n \in [N]} a_n e(n\theta_i)|^2 \leq (N + \delta^{-1})\sum_{n \in [N]} |a_n|^2.\end{equation}
\end{theorem}
As it stands, this inequality may be thought of as a kind of approximate Bessel's inequality; the separation condition \eqref{separation} acts to ensure that the exponentials $e(n\theta_i)$, $i = 1,\dots,k$, are roughly orthogonal over $n \in [N]$. The rather clean form of \eqref{l-sieve} was proved by Montgomery and Vaughan \cite{montgomery-vaughan} and to prove it somewhat careful arguments are needed. It is much easier to establish a weaker inequality in which the right-hand side is replaced by (say) $8(N + \delta^{-1})\sum |a_n|^2$; at the level of the discussions in this paper, this is just as good.

The large sieve gets its name from the fact that \eqref{l-sieve} may be used to give bounds for certain sieving problems, of which the following is an example.

\begin{problem}[Sieve problem]\label{sieve-problem} Let $N \geq 1$. Suppose that for each off prime $p \leq \sqrt{N}$ one is given a set $S_p \subseteq \Z/p\Z$ with $|S_p| = (p-1)/2$. Let $A \subseteq \{1,\dots,N\}$ be the set obtained by ``sieving out'' the residue classes $S_p$ for each prime, that is to say by removing from $[N]$ each $n$ for which  $n \md{p}$ lies in $S_p$ for some $p$. What upper bounds can one place on $A$?\end{problem}

\emph{Remark.} In actual fact one may consider far more general settings. The \emph{sieving limit} $X = \sqrt{N}$ may be reduced, one might sieve only by a subset of the primes $p \leq X$ rather than by all primes, and the size of $S_p$ might vary less regularly with $p$. 

The large sieve gives, by an argument of Montgomery \cite{montgomery-lsieve}, the following result concerning Problem \ref{sieve-problem}. 

\begin{theorem}[Large sieve bound]\label{l-sieve-bound}
In Problem \ref{sieve-problem} we have the bound $|A| \leq C\sqrt{N}$ for some absolute constant $C$.\endproof
\end{theorem}

Theorem \ref{l-sieve-bound} is essentially sharp. Indeed if one takes $S_p$ to consist of the quadratic non-residues for each prime $p$ then the unsieved set $A$ contains all of the squares less than or equal to $N$, and so $|A| \geq (1 -o(1))\sqrt{N}$. It is, however, very hard to think of an essentially different example giving a comparable lower bound. As a result Helfgott and Venkatesh \cite{helfgott-venkatesh} and independently Croot and Elsholtz (cf. \cite[Problem 7.4]{croot-lev}) were motivated to make the beautiful guess that any set $A$ of size close to $\sqrt{N}$ which survives the sieving process in Problem \ref{sieve-problem} is essentially the set of values of some quadratic polynomial.

\begin{conjecture}[Inverse conjecture for the large sieve]\label{sieve-inverse}
Suppose that the set $A$ is the result of the sieving process in Problem \ref{sieve-problem} and that $|A| \geq N^{0.499}$. Then all but $O(N^{o(1)})$ points of $A$ are contained in the set of values of a quadratic polynomial $f(n) = an^2 + bn + c$.
\end{conjecture}

\emph{Remark.} In fact it is not unreasonable to think, as Croot, Elsholtz, Helfgott and Venkatesh did, that much more should be true. One might consider more general sieving situations than the one in Problem \ref{sieve-problem}, and even under an assumption as weak as $|A| \geq N^{\eta}$ the set $A$ should have some fairly rigid ``algebraic'' structure.

Conjecture \ref{sieve-inverse} seems to be of interest in its own right. Furthermore suitable variants of it ought to have applications, for example to Ostmann's Inverse Goldbach Problem as considered in a series of papers by Elsholtz \cite{elsholtz}.

We are not able to establish Conjecture \ref{sieve-inverse} or anything close to it. The aim of this paper is to develop ideas which go somewhat beyond the large sieve in using the \emph{structure} of the residue classes $S_p$ rather than merely their size. Rather than develop these ideas in the most general context we use them to address the most extreme case in which each set $S_p$ of residues is an interval in $\Z/p\Z$. The following is our main result.

\begin{theorem}[Interval sieve]\label{middle-half}
Let $N \geq 1$, and suppose that for each odd prime $p \leq \sqrt{N}$ one is given an interval $I_p \subseteq \Z/p\Z$ of length $(p-1)/2$. Let $A \subseteq \{1,\dots,N\}$ be the set obtained by sieving out all $n$ for which $n\mdlem{p}$ lies inside $I_p$ for some $p$. Then $|A| \ll N^{1/3 + o(1)}$.
\end{theorem}
\emph{Remarks.} We have concocted this problem so that the theorem is not obviously trivial given known results. For the more specific situation in which, for example, $I_p = [p/4,3p/4]$ is the ``middle-half'' it follows from results of Jutila \cite{jutila} that $|A| = O(1)$, and it is quite possible that this bound could be made rather effective using the results of Granville and Ramar\'e \cite{granville-ramare}. It seems quite reasonable to suppose that $|A| = O(1)$ in the somewhat more general setting of Theorem \ref{middle-half}, uniformly in the choice of the intervals $I_p$.  

Our own arguments are fairly routine, but we do import a very interesting result of Roberts and Sargos \cite{robert-sargos} concerning the spacing properties of the set of unit fractions $\{1/x : X \leq x < 2X\}$. If the argument were written down in a self-contained manner, this result would certainly be the beef. As we shall sketch later on, the \emph{exponent pairs conjecture} may be used to obtain the stronger bound $|A| \ll N^{o(1)}$.

\section{A variant of the large sieve}

The goal of this section is to prove the following proposition, which is our variant of the large sieve. In this proposition $\mathbf{a} = (a_n)_{n \in [N]}$ is a sequence of complex numbers and we write
\[ \Vert \mathbf{a} \Vert_r := \big( \sum_{n \in [N]} |a_n|^r \big)^{1/r}\] for $r \geq 1$.

\begin{proposition}\label{prop3.1}
Let $\mathbf{a} = (a_n)_{n \in [N]}$ be a sequence of complex numbers. Let $X,N$, $1 \leq X \leq N$, be parameters. Then 
\[ \sum_{X \leq x < 2X} |\sum_{n \in [N]} a_n e(n/x)|^2 \ll (N + X^3)^{1/4}X^{1/2 + o(1)} \Vert \mathbf{a} \Vert_{4/3} \Vert \mathbf{a} \Vert_1.\]
The same is true with $e(n/x)$ replaced by $e(2n/x)$.
\end{proposition}
\emph{Remark.} When we apply this proposition to the interval sieve problem in the next section we will need the exponentials $e(2n/x)$ as well and that is why we have mentioned them here. One could use $e(kn/x)$ for any fixed $k$, though the implied constants would depend on $k$.

\proof The proof of the corresponding inequality when $e(2n/x)$ replaces $e(n/x)$ is identical and we say nothing more about it. Let $\psi : \R \rightarrow \R_{\geq 0}$ be a Beurling-Selberg function with the following properties:

\begin{enumerate}
\item $\psi(t) \geq 1$ for $|t| \leq N$;
\item $\tilde \psi(\xi)$ is supported on $|\xi| \leq 1/X^3$, where $\tilde\psi(\xi) := \int^{\infty}_{-\infty} \psi(x) e^{-i\xi x}\, dx$;
\item $\Vert \psi \Vert_1 \ll N + X^3$.
\end{enumerate}

The use of these majorants in analytic number theory is well-known, and the book of Montgomery \cite{montgomery-ten-lecture} or the article of Vaaler \cite{vaaler} may be consulted for more information concerning their construction.

Now the left-hand side in our proposition may be expanded as
\[ \sum_{n \in [N]} \frac{a_n}{\psi(n)^{1/4}} \psi(n)^{1/4} \sum_{m,x} \overline{a_m} e(n/x)e(-m/x).   \]
By H\"older's inequality this is at most
\[ \big(\sum_{n \in [N]} \frac{|a_n|^{4/3}}{\psi(n)^{1/3}}\big)^{3/4} \big( \sum_n \psi(n)|\sum_{m,x} \overline{a_m}  e(n/x) e(-m/x)|^4 \big)^{1/4}.\]
By property (i) of the majorant $\psi$, the first factor is bounded by $\Vert \mathbf{a}\Vert_{4/3}$.
The expression inside the second bracket may be expanded as
\begin{align*} \sum_{x_1,x_2,x_3,x_4} \sum_{m_1,m_2,m_3,m_4} \overline{a_{m_1}a_{m_2}}a_{m_3}a_{m_4} & e(-m_1/x_1) e(-m_2/x_2) e(m_3/x_3)e(m_4/x_4)\times \\ & \times\sum_n \psi(n) e(n (\frac{1}{x_1} + \frac{1}{x_2} - \frac{1}{x_3} - \frac{1}{x_4})),\end{align*} which is bounded by
\begin{equation}\label{bounded-by} \Vert \mathbf{a} \Vert_1^4\sum_{x_1,x_2,x_3,x_4}  |\widehat{\psi} (\frac{1}{x_1} + \frac{1}{x_2} - \frac{1}{x_3} - \frac{1}{x_4})|.\end{equation} Here, the hat denotes the Fourier transform on $\Z$, so we are writing
\[ \widehat{\psi}(\theta) := \sum_n \psi(n) e^{2\pi i n \theta}.\]
By the Poisson summation formula we have 
\[ \widehat{\psi}(\theta) = \sum_n \tilde \psi(\theta - n),\]
and so by properties (ii) and (iii) of $\psi$ we see that \eqref{bounded-by} is at most $C\Vert a \Vert_1^4 (N + X^3)$ times the number of quadruples $x_1,x_2,x_3,x_4 \in [X,2X)$ with 
\[ |\frac{1}{x_1} + \frac{1}{x_2} - \frac{1}{x_3} - \frac{1}{x_4}| \leq \frac{1}{X^3}.\]
It follows from Theorem 2 of Roberts and Sargos \cite{robert-sargos} that there are $\ll X^{2 + o(1)}$ such quadruples. The proposition follows quickly.\endproof

\emph{Remark.} The reader familar with basic duality theory in Banach spaces may recognise some aspects of the proof of Proposition \ref{prop3.1}. It is very closely modelled on the proof that $\Vert T \Vert_{4/3 \rightarrow 2} = \Vert T^* \Vert_{2 \rightarrow 4}$, where $T : B(X) \rightarrow B(Y)$ and $T^* : B(Y) \rightarrow B(X)$ are mutually adjoint operators on spaces $B(X),B(Y)$ of bounded functions. One might also regard our variant of the large sieve as a kind of \emph{restriction theorem}, where one looks at the Fourier transform of the sequence $(a_n)$ restricted to the set of frequencies $\{1/x : x \in [X,2X)\}$. Proposition \ref{prop3.1} then reflects a kind of ``discrete curvature'' of this set of frequencies, and is closely analogous to such estimates as the Tomas-Stein restriction theorem (cf. \cite{tao-restriction}).

\section{Interval sieve problem}

It is a reasonably straightforward matter to apply Proposition \ref{prop3.1} to get the stated bound for the interval sieve problem.

\begin{lemma}\label{large-fourier} Suppose that $p$ is a prime, that $I_p \subseteq \Z/p\Z$ is an interval of length $(p-1)/2$ and that $A \subseteq [N]$ is a set such no $a \in A$ has $a\mdlem{p} \in I_p$. Then either $\sum_{n \in [N]} 1_A(n) e(n/p)$ or $\sum_{n \in [N]} 1_A(n) e(2n/p)$ has magnitude at least $|A|/3$.
\end{lemma}
\proof We note that $1 - 2\cos \theta + \cos 2\theta  \leq 0$ when $|\theta| \leq \pi/2$; rewriting the left-hand side as $2\cos\theta (\cos \theta - 1)$, this becomes clear. It follows that there is $\beta \in [0,1]$ (depending on $I_p$) such that if $n \in A$ then
\[ 1 - 2\cos 2\pi (\frac{n}{p} + \beta) + \cos 4\pi (\frac{n}{p} + \beta) \leq 0,\] and hence
\[ 1 - e(\beta)e(n/p) - e(-\beta)e(-n/p) + \frac{1}{2}e(2\beta)e(2n/p) + \frac{1}{2}e(-2\beta)e(-2n/p) \leq 0.\]
Summing over $n \in A$ and using the triangle inequality, one obtains
\[ |A| \leq 2|\sum_{n \in [N]} 1_A(n) e(n/p)| + |\sum_{n \in [N]} 1_A(n) e(2n/p)|,\] from which the result follows immediately.\endproof

\emph{Remark.} What we have shown here is that if $A \md{p}$ does not meet $I_p$ then $A$ has an extremely large discrete Fourier coefficient in $\Z/p\Z$. In the usual application of the large sieve, one shows that the $L^2$-mass of the discrete Fourier transform of $A$ in $\Z/p\Z$ has significant mass away from the zero mode; this is much weaker information and in our case is wasteful as it does not utilise the specific additive structure of the excluded residues $I_p$.

We may now prove Theorem \ref{middle-half}. Suppose that $A$ is the set of those elements of $[N]$ which remain after sieving by all residues in $I_p$, for all $p \leq \sqrt{N}$. Then, setting $a_n := 1_A(n)$, the previous lemma implies that
\[ \max ( |\sum_{n \in [N]} 1_A(n) e(n/p)|, |\sum_{n \in [N]} 1_A(n) e(2n/p)| ) \geq |A|/3.\]
Note that if $\mathbf{a} = (a_n)_{n \in [N]}$ then $\Vert \mathbf{a} \Vert_1 = |A|$ and $\Vert \mathbf{a} \Vert_{4/3} = |A|^{3/4}$. Substituting into Proposition \ref{prop3.1}, we obtain for any $X$ the bound
\[ X |A|^2 \ll (N + X^3)^{1/4} X^{1/2 + o(1)}|A|^{7/4}.\] Taking $X = N^{1/3}$ leads to the stated bound.\endproof

\emph{Remark.} Our argument has something in common with the argument used to obtain lower bounds in the Kakeya problem from restriction estimates for the sphere, which had its origin in the work of Fefferman \cite{fefferman}. Indeed our sieve bound, phrased differently, provides a \emph{lower bound} of $N - O(N^{1/3 + \eps})$ for a union of sets $(I_p + p\Z) \cap \{1,\dots,N\}$ which, at a stretch, might be thought of as ``lines'' in different directions in the spirit of the Kakeya problem. These ideas have been considered in a number theoretic context before in the work of Bourgain \cite{bourgain-montgomery}.

\section{Further remarks}

We conjecture that the following is true.

\begin{conjecture}[Spacing of unit fractions]\label{conj-5} Suppose that $r \geq 1$ is an integer. Then the number of $x_1,\dots,x_{2r} \in [X,2X)$ such that
\begin{equation}\label{eq56} |\frac{1}{x_1} + \dots + \frac{1}{x_r} - \frac{1}{x_{r+1}} - \dots - \frac{1}{x_{2r}}| \leq \frac{1}{X^{r+1}}\end{equation} is $\ll_{r,\eps} X^{r + \eps}$ for all $\eps > 0$.
\end{conjecture}

If this did hold for a particular value of $r$ then a straightforward modification of our arguments (using exponents $p = 2r/(2r-1)$ and $q = 2r$ in H\"older's inequality in the proof of Proposition \ref{prop3.1}) would lead to a bound $|A| \ll  N^{\frac{1}{r+1} + o_r(1)}$ in Theorem \ref{middle-half}. By bounding the number of solutions to \eqref{eq56} using a $2r$-power moment of exponential sums as in \cite{robert-sargos} one may confirm that the conjecture would follow if we had a bound
\begin{equation}\label{conjectured-bound} |\sum_{X \leq x < 2X} e(\xi/x)| \ll X^{1/2 + o_r(1)}\end{equation} for $X^2 \leq |\xi| \leq X^{r+1}$.
Such a bound is a consequence of the so-called \emph{exponent pairs hypothesis,} stated on p.214 of \cite{iwaniec-kowalski}. An excellent source of information on exponent pairs is the book \cite{graham-kolesnik}, though the exponent pairs hypothesis itself is conspicuously absent from that book. Nonetheless it seems to be fairly widely believed, and in any case the bound \eqref{conjectured-bound} accords with the commonly-held belief that exponential sums should exhibit square-root cancellation unless there is a ``good'' reason for them not to.  

Improving the bounds for $\sup_{\xi} |S(\xi)|$ seems closely related to the Lindel\"of hypothesis, which is equivalent to proving that 
\[ |\sum_{X \leq x < 2X} e(\xi \log x)| \ll X^{1/2 + o_r(1)}\] for $|\xi| \ll X^r$, for all $r$.
We do not intend this remark to be taken too seriously: it stems from the observation that the derivatives of the phase $\log x$ are the same as those of $1/x$, and the derivative structure of a phase is often important in the estimation of exponential sums. Nevertheless, it would be very surprising if anything close to \eqref{conjectured-bound} were proved tomorrow. Perhaps better evidence for this is that an affirmative solution to the Dirichlet divisor problem and the Gauss circle problem would follow from \eqref{conjectured-bound} for $r = 3$. See \cite[Chapter 4]{graham-kolesnik} for more information.

We note that \eqref{conjectured-bound} is not known for $r = 2$; Roberts and Sargos bound the 4th moment of their exponential sums directly, without obtaining a bound for the supremum of those sums.

\section{Acknowledgements}

I would like to thank Roger Heath-Brown for drawing my attention to the reference \cite{robert-sargos} and Christian Elsholtz for helpful remarks.

\providecommand{\bysame}{\leavevmode\hbox to3em{\hrulefill}\thinspace}
\providecommand{\MR}{\relax\ifhmode\unskip\space\fi MR }
\providecommand{\MRhref}[2]{%
  \href{http://www.ams.org/mathscinet-getitem?mr=#1}{#2}
}
\providecommand{\href}[2]{#2}

     \end{document}